
\magnification\magstep1
\scrollmode                   


\font\eightrm=cmr8            
\font\eightit=cmti8           


\def\cite#1{{\rm[#1]}}        
\def\opname#1{\mathop{\rm#1}\nolimits} 
\def\refno#1. #2\par{\smallskip\item{[#1]} #2\par}



\newtoks\rightheadtext      
\newtoks\leftheadtext       

\rightheadtext{A TANGENT BUNDLE ON DIFFEOLOGICAL SPACES}

\leftheadtext{CARLOS A. TORRE}

\headline={\ifnum\pageno>1  
            \ifodd\pageno   
             \hfil \eightrm \the\rightheadtext \hfil
             \llap{\tenrm\folio}%
            \else \rlap{\tenrm\folio}
             \hfil \eightrm \the\leftheadtext \hfil \fi
           \else\hfil \fi}
\footline={\hfil}           


\outer\def\section#1. #2\par{
      \bigskip\bigskip       
      \message{#1. #2}
      \leftline{\bf#1. #2}
      \nobreak\smallskip     
      \noindent}             

\outer\def\varsection#1\par{
      \bigskip\bigskip       
      \message{#1 }
      \leftline{\bf#1 }
      \nobreak\smallskip     
      \noindent}             

\def\declare#1. #2\par{
      \medskip\noindent     
      {\bf#1.}\rm           
      \enspace\ignorespaces 
      #2\par\smallskip}     

\def\demo#1. {\smallskip    
      \noindent             
      {\bf #1}.
      \enspace\ignorespaces}
\def\enddemo{\par\smallskip} 



\def\al{{\cal G}}             
\def\D{{\cal D}}              
\def\M{{\cal M}}              
\def\Tn{{\cal T}}             


\def\N{{\bf N}}               
\def\R{{\bf R}}               


\def\Alt{\opname{Alt}}        
\def\Diff{\opname{Diff}}      
\def\Hom{\opname{Hom}}        
\def\im{\opname{im}}          


\def\dR{{\rm dR}}             


\def\SDS{{\eightrm SDS}}      


\def\ox{\otimes}              
\def\W{\Lambda}               
\def\w{\wedge}                
\def\x{\times}                
\def\.{\cdot}                 
\def\:{\colon}                



\def\ca{C^\infty}             
\def\cader{{\cal C}^\infty}   
\def\oa{\Omega_{F_0}}         
\def\pp{\tilde p}             
\def\sptilde{\widetilde{\;}\,}


\def\qed{\allowbreak\qquad\null
         \nobreak\hfill\square} 
\def\square{\vrule height 5pt depth 0pt
             width 5pt}       



\def\sepword#1{\qquad\hbox{#1}\quad} 


\def\set#1{\{\,#1\,\}}        
\def\<#1>{\langle#1\rangle}   


\def\frac#1#2{{#1\over#2}}    


\def\Adams{1}            
\def\DonatoGeometrie{2}
\def\DonatoVirasoro{3}
\def\DonatoFibrations{4} 
\def\KirillovBook{5}
\def\SouriauAlgoritme{6}
\def\Woodhouse{7}


\topglue 1.5cm

\centerline{\bf A Tangent Bundle on Diffeological Spaces}

\bigskip

\centerline{Carlos A. Torre}
\bigskip
\centerline{\it Escuela de Matem\'atica,
                Universidad de Costa Rica,}
\centerline{\it 2060 San Jos\'e, Costa Rica}
\bigskip

{\narrower
\eightrm
\baselineskip=9.5pt
\noindent{\eightit ABSTRACT}.
We define a subcategory of the category of diffeological spaces, which
contains smooth manifolds, the diffeomorphism subgroups and its
coadjoint orbits. In these spaces we construct a tangent bundle,
vector fields and a de~Rham cohomology.
\par
\smallskip}

\vfootnote{}{{\it Keywords\/}: Smooth diffeology, tangent space,
diffeomorphism group.}
\vfootnote{}{1991 {\it Mathematics Subject Classification\/}: 57R55,
58B10; 22E65.}

\bigskip


\section 1. Introduction

The category of diffeological spaces
\cite{\DonatoGeometrie,~\SouriauAlgoritme} extends that of manifolds
and allows many topological and geometrical constructions, including
products, quotients, forms, homotopies,
fibrations~\cite{\DonatoFibrations}. It has the property that
$\Hom(X,Y)$ is an object in a canonical way whenever $X$ and $Y$ are,
thereby allowing the study of infinite dimensional objects.

It is known that symplectic manifolds play a central role in the
modelling of physical systems~\cite{\Adams,~\Woodhouse}, and in
particular coadjoint orbits have a natural symplectic
structure~\cite{\KirillovBook}. By using a covariant definition of
forms, this result is extended by Souriau also to the category of
diffeological spaces and even allows the construction of a
prequantization~\cite{\DonatoVirasoro,~\SouriauAlgoritme} whenever
a certain cohomological obstruction is zero.

In this paper we introduce the concept of a smooth diffeological
space~(\SDS). This is a subcategory that includes manifolds, but is
general enough to include groups of diffeomorphisms and
infinite-dimensional coadjoint orbits. This category allows a
dynamical modelling: we prove that these objects possess canonical
tangent bundles and that each of them are also \SDS. This allows the
construction of flows and Lie algebras of vector fields, and of
de~Rham cohomologies.

\section 2. Smooth diffeologies

Given a set $X$, an {\it $n$-plaque\/} of~$X$ is a function
$p\: U \to X$ where $U$ is an open subset of~$\R^n$ containing the
origin~$0$. A {\it diffeology\/} on $X$~\cite{\SouriauAlgoritme} is a
set $P(X)$ of $n$-plaques for every~$n$ such that:
\itemitem{(i)}
  The images of the plaques cover $X$.
\itemitem{(ii)}
  If a set $(p_i)$ of $n$-plaques admits a common extension,
  then the smallest such extension is also an $n$-plaque in~$P(X)$.
\itemitem{(iii)}
  For every $\psi \in \ca(U',U)$ where $U$, $U'$ are open subsets of
  $\R^m$, $\R^n$ respectively, and every plaque $p\: U \to X$,
  $p \circ \psi$ is also in~$P(X)$.

\smallskip

The set of $n$-plaques is denoted $P^n(X)$, and the set of $n$-plaques
$p$ such that $p(0) = F \in X$ is denoted $P_F^n(X)$. It may be seen
that every set $S$ of plaques on $X$ generates a smallest diffeology
$P_S$ which is formed by plaques $p\: U \to X$ such that for all
$r \in U$ there exists $U_r \subset U$ open with
$p|_{U_r}= f \circ \phi$ where $\phi \: U_r \to V \subset \R^n$ is
smooth and $f \in S$.

A map $f \: X \to Y$ between two spaces with diffeologies $P(X)$
and $P(Y)$ is called {\it differentiable\/} if $p \in P(X)$
implies $f \circ p \in P(Y)$. The set of such maps is denoted by
$\D^\infty(X,Y)$.

We shall now introduce a particular class of diffeologies that we call
{\it smooth}, in order to set the stage for infinite-dimensional
tangent spaces in the diffeological category.

\declare Definition 2.1.
Let $k \in \N$, a triple $(X, P(X), \sim)$ is called a
{\it $C^k$-diffeology\/} (or {\it smooth diffeology\/} when
$k = \infty$) if:
\itemitem{(a)}
The pair $(X, P(X))$ is a diffeology.
\itemitem{(b)}
$\sim$ is a collection $\set{\sim_F^n : 1 \leq n \leq k,\ F \in X}$,
where $\sim_F^n$ is an equivalence relation on the set $P_F^n(X)$,
that satisfy a consistency condition:
$p_1 \circ \phi \sim_F^m p_2 \circ \phi$ whenever $p_1 \sim_F^n p_2$
and $\phi \in C^\infty(U',U)$, with $U' \subset \R^m$ and
$p_1,p_2\: U \to X$; and moreover, $p_1 \sim_F p_1|_V$ whenever
$p_1 \: U \to X$ and $V \subset U$ is an open neighborhood of~$0$.

\itemitem{(c)}
the smooth diffeology is called {\it linear\/} if the set
$V := \set{V_F^n : F \in X,\ n = 1,\dots,k}$ where
$V_F^n = P_F^n(X) /{\sim_F^n}$, carries a vector space structure that
satisfies the following consistency condition: whenever
$p_{12} \in [p_1] + [p_2]$ with $p_i\: U \subset \R^n \to X$ and
$\phi \: U' \subset \R^m \to U$ then
$p_{12} \circ \phi \in [p_1 \circ \phi] + [p_2 \circ \phi]$.

The class of the plaque $p(t)$ at a point $F \in X$ will be denoted by
$[p]$ or by $[p]_t$. Whenever there are several spaces, we may use the
notation $[p]^X$ instead.

The linear structure is called {\it continuous\/} if given two
$(n+m)$-plaques $p_i(r,s)$ with $i = 1,2$ such that
$p_1(r,0) = p_2(r,0)$ then there exists a plaque $p_{12}(r,s)$ such that
$$
[p_{12}(r,s)]_s = [p_1(r,s)]_s + [p_2(r,s)]_s  \qquad\hbox{for all $r$}.
$$

\declare Definition 2.2.
The set $\Tn_F^n X := P_F^n(X) /{\sim_F^n}$ is called the $n$-th
{\it tangent space\/} at~$F$, and the disjoint union
$\Tn^n X := \bigsqcup_{F\in X} \Tn_F^n X$ is called the $n$-th
{\it tangent bundle\/} over~$X$.

Notice that $\Tn_F^n X$ need not carry a linear structure. For example,
the union of two smooth curves that intersect transversally at~$F$ is
a smooth diffeological space, for which $\Tn_F^n X$ is the union of two
lines.

\declare Definition 2.3.
Let $(X,P(X),\sim)$, $(Y,P(Y),\approx)$ be smooth diffeologies; a
differentiable function $f \: X \to Y$ is called {\it smooth\/} at~$F$
if for all $n \in \N$, $f \circ p_1 \approx_{f(F)}^n f \circ p_2$
whenever $p_1 \sim_F^n p_2$.

If the smooth diffeologies are linear, $f$ is called a {\it smooth
map\/} if, in addition,
$$
D_F^n f : \Tn_F^n X \to \Tn_{f(F)}^n Y : [p] \mapsto [f \circ p]
$$
is linear for each $F \in X$. The set of smooth functions is denoted
$\cader(X,Y)$. Notice that $\cader(X,Y) \subset \D^\infty(X,Y)$. In
particular, if $Y = \R$ with the smooth manifold diffeology that is
described below, then the notation $\cader(X)$ will be used instead.

Notice that the set of smooth diffeological spaces as objects with the
smooth maps as morphisms is a category, since if $f \: A \to B$ and
$g \: B \to C$ are smooth then $g \circ f$ is differentiable and
preserves the equivalence relations. Moreover, the chain rule
$$
D^n_F(g \circ f)[p] = D^n_{f(F)}g (D^n_F f[p])
$$
holds, and linearity is preserved by this composition; therefore
$g \circ f$ is a smooth map if $f$ and $g$ are smooth maps.

Let $(X,P,\sim)$ be a linear \SDS\ (smooth diffeological space) and
let $Y \subset X$. The {\it subspace diffeology\/} is formed by the
plaques $p$ such that $p(r) \in Y$ for all~$r$. Restrict $\sim$ to
$P(Y)$. $Y$ is called a sub-\SDS\ if $P^n(Y)/{\sim_F^n}$ is a subspace
of $P^n(X)/{\sim_F^n}$ for each $F \in Y$.

Let $(X,P(X),\sim)$ and $(Y,P(Y),\equiv)$ be two (linear) \SDS. Denote
by $P(X \x Y)$ the product diffeology on $X \x Y$ formed by plaques of
the form $p(r) = (p_1(r), p_2(r))$ where $p_1 \in P(X)$ and
$p_2 \in P(Y)$. Define
$$
(p_1,p_2) \approx_F^n (p'_1,p'_2)  \sepword{if and only if}
 p_1 \sim_{F_1}^n p'_1 \quad\hbox{and}\quad p_2 \equiv_{F_2}^n p'_2.
$$
Then $V^n_{(F_1,F_2)} = V^n_{F_1} \x V^n_{F_2}$. The triple
$(X \x Y, P(X \x Y), \approx)$ is an \SDS.

\medskip

Some examples of \SDS\ are the following:
\smallskip

\item{(1)}
Let $M$ be a smooth manifold modelled on a locally convex vector
space~$V$. Let $P(M)$ be the manifold diffeology (formed by the smooth
plaques). For each $F \in M$, choose a chart $(U_F,\alpha_F)$
around~$F$, and define $p_1 \sim_F^n p_2$ if
$$
D^m(\alpha_F \circ p_1)(0) = D^m(\alpha_F \circ p_2)(0)
 \sepword{for all}  m \leq n.
$$
The map
$$
\Tn_FM \to V : [p] \mapsto \frac d{dt}(\alpha_F\circ p)(t)\biggr|_{t=0}
$$
is a bijection and defines $V_F^1$. Other spaces are defined
similarly. In this case $\ca(M) \subset \cader(M) = \D^\infty(M)$,
where $\ca(M)$ denotes the set of smooth functions with respect to
the manifold structure. In particular, when $M = \R^n$, we recover the
standard differentiable structure on a finite-dimensional vector space.

\item{(2)}
Let $M$ be a smooth finite dimensional manifold and let
$X = \Diff_c(M)$ (the infinite-dimensional group of diffeomorphisms of
$M$ with compact support). Define $P(X)$ as the set of functions
$p \: U \to X$, $U \subset \R^n$ open, such that
$$
\phi_p : U \x M \to M : (r,m) \mapsto p(r)(m)
$$
is smooth. Given $g \in X$ and
$p_1,p_2\in P^n(X)$ with $p_1(0) = p_2(0) = g$, we define
$p_1 \sim_g^n p_2$ if
$$
D^i p_1(t)(m)\bigr|_{t=0} = D^i p_2(t)(m)\bigr|_{t=0}
  \sepword{for all}  m \in M, \ i \leq n.
$$
Then $\Tn_g^1 X = \Gamma_c(M)$, the space of vector fields with
compact support on~$M$.

\item{(3)}
There is a natural smooth diffeology on coadjoint orbits. Let $G$ be a
subgroup of diffeomorphisms (of the group $X$ defined above) with Lie
algebra $\al$ such that $\frac d{dt} p(t)|_{t=0} \in \al$ for each
diffeotopy~$p$ in~$G$. This is a condition that holds for example if
$G$ is the group $\Diff_c(M)$ or the group of symplectic
diffeomorphisms, or if $G$ is finite dimensional (it is an open
question whether this holds for any closed subgroup). On any coadjoint
orbit $\oa$ define the diffeology $P(\oa)$ generated by the set of
plaques of the form $b(r) = K(p(r))F$ where $F \in \oa$ and
$p \in P(G)$ (defined above from the subspace diffeology) and $K$ is
coadjoint action. In $P(\oa)$ define $b_1 \sim_F^n b_2$ if
$$
D^m b_1(t)(Y)\bigr|_{t=0} = D^m b_2(t)(Y)\bigr|_{t=0}
$$
for each $Y$ in $\al$, and for each $m \leq n$. Let $b \in P(\oa)$,
and let
$$
b(t) =: K(p(t))F
$$
where $p \in P(G)$. Define $\xi = (d/dt) p(t)|_{t=0}$. Then the map
$[b] \mapsto dK(\xi)F$ is a bijection between $\Tn_F\oa$ and
$\al/\al(F)$, and allows us to regard $\Tn_F\oa$ as a vector space by
transport of structure. In this way $\oa$ becomes a $C^1$ linear
diffeological space.

We shall prove that $\Tn^m X$ has a natural smooth diffeology.
Consider the set of maps of the form
$$
\pp(r_1,\dots,r_n) := [p(r_1,\dots,r_n,s_1,\dots,s_m)]_s
$$
where $p$ is a $(n+m)$-plaque. Let $P(\Tn^m X)$ be the set of plaques
generated by these plaques (via condition~(ii) of the definition of
diffeology).

Any $[\alpha] \in \Tn^mX$ is a class of $m$-plaques at some point~$F$.
Let $\alpha\: U \to X$, $U\subset \R^m$, $\alpha(0) = F$ be such a
plaque, and let $\pp_1$, $\pp_2$ be two $n$-plaques on $\Tn^mX$ such
that $\pp_1(0) = \pp_2(0) = \alpha$; define
$$
\pp_1 \approx_{[\alpha]}^n \pp_2  \iff  p_1 \sim_F^{n+m} p_2.
$$

\proclaim Proposition 2.1.
$(\Tn^mX, P, \approx)$ is a smooth diffeology.

\demo Proof.
These plaques cover $\Tn^m X$ since, in particular, if
$[p(s_1,\dots,s_m)] \in \Tn^m X$, then $\tilde q(0) = [p]$, where
$q(t,s_1,\dots,s_m) := p(s_1,\dots,s_m)$.

Now let $\psi \: U' \to U$ where $U' \subset \R^k$, $U \subset \R^n$
are open. Then
$$
\eqalign{
\pp(\psi(r_1,\dots,r_k))
&= [p(\psi(r_1,\dots,r_k),s_1,\dots,s_m)]_s
\cr
&= [(p \circ (\psi \ox 1_m))\sptilde(r_1,\dots,r_k,s_1,\dots,s_m)]_s,
\cr}
$$
where the map
$$
(\psi \ox 1_m)(r_1,\dots,r_k, s_1,\dots,s_m)
 := (\psi(r_1,\dots,r_k), s_1,\dots,s_m)
$$
is smooth, $p \circ (\psi \ox 1_m)$ is a $(k + m)$-plaque and
$\pp\circ \psi$ is a $k$-plaque. Therefore $P(\Tn^mX)$ is a diffeology.

Let $\psi \: U' \to U$ and assume $\pp_1 \approx_{[\alpha]}^n \pp_2$,
then $p_1 \sim_F^{n+m} p_2$, therefore
$$
p_1 \circ (\psi \ox 1_m) \sim_F^{n+m} p_2 \circ (\psi \ox 1_m),
$$
thus
$$
p_1 \circ (\psi \ox 1_m)\sptilde \approx_{[\alpha]}^n
 p_2 \circ (\psi \ox 1_m)\sptilde.
$$
It follows that
$\pp_1 \circ \psi \approx_{[\alpha]}^n \pp_2 \circ \psi$.   \qed
\enddemo

The space $\Tn_F X$ is given the subspace diffeology of $\Tn X$. For
example, for $n = 1$ it consists of plaques $\pp$ of the form
$\pp\: U \to \Tn_F X$ such that $\pp(r) = [p(r,t)]_t^X$ where
$p(r,0) = F$ for each $r$; then $p(r,t)$ is a $t$-curve through $F$
for each~$r$.

If $X$ has a linear structure, then we can define a linear structure on
$\Tn^n X$ in the following way: take a point $(F,[\alpha])$ at
$\Tn^n X$ and two $k$-vectors at this point (these are classes of
$k$-plaques $[\pp_i] \: U \subset \R^k \to \Tn^n X$ such that
$s \mapsto p_i(0,s) \in [\alpha]$. Then define
$$
[\pp_1] + c[\pp_2] = [\pp_{12}],
$$
where $p_{12} \in [p_1] + c[p_2] - c[\bar\alpha]$, where
$\bar\alpha(r,s) = \alpha$ for all~$r$ (notice that $[\pp_{12}]$ is a
$k$-plaque through $(F,[\alpha])$ and that if $[\pp_1] = 0$ then
$\pp_1 = \tilde\alpha$). One can check that $\Tn_F^n X$ becomes
also a linear~\SDS.

\proclaim Proposition 2.2.
Let $f \: X \to Y$ be smooth; then $D^mf \: \Tn^m X \to \Tn^m Y$ is
smooth for each~$m$. Also,
$D^m_{F_0} \: \Tn^m_{F_0}X \to \Tn^m_{f(F_0)} Y$ is smooth for each
$m$ and each $F_0\in X$.

\demo Proof.
First, we shall prove that $D^m f$ is differentiable. Let $\pp$ be an
$n$-plaque on $\Tn^m X$, then $\pp(r) = [p(r,s)]_s^X$, therefore
$f \circ \pp$ is an $(n+m)$-plaque on $Y$ and
$[f\circ p(r,s)]^Y_s \in \Tn^mY$ for all~$r$, and it is equal to
$(f \circ p)\sptilde (r)$, which is an $n$-plaque on $\Tn^m Y$.
Therefore
$$
(D^mf) \circ \pp(r) = [f \circ p(r,s)]_s^Y = (f \circ p)\sptilde (r).
$$
This proves that $D^mf$ is differentiable. It is also smooth: if
$\pp_1 \approx \pp_2$ at $(F,[\alpha])$ then $p_1 \sim p_2$, therefore
$p_1 \sim f \circ p_2$ and
$(f \circ p_1)\sptilde \approx (f \circ p_2)\sptilde$ at
$(f(F), [f \circ \alpha])$. Thus
$D^m f \circ \pp_1 \approx D^m f\circ \pp_2$ at
$(f(F), [f \circ \alpha])$.

If $X$ has a linear structure, then one can check that the
differential of these maps are also linear.   \qed
\enddemo

\proclaim Proposition 2.3.
The projection $\pi \: \Tn X \to X$ defined by $\pi[\alpha]_F := F$ is
smooth.

\noindent{\it Proof}.\enspace
Let $\pp$ be a $k$-plaque of $\Tn X$ on $(F,[\alpha])$,
$\pp(r) = [p(r,t)]_t$. Then $(\pi \circ \pp)(r) = p(r,0)$ is a plaque
on~$X$. Therefore $\pi$ is differentiable. Assume
$\pp_1,\pp_2 \in P(\Tn X)$. If $\pp_1 \approx \pp_2$ then
$p_1(r,t) \sim^{n+1} p_2(r,t)$, then $p_1(r,0) \sim^n p_2(r,0)$;
therefore $\pi (\pp_1) \sim \pi (\pp_2)$. If $X$ has a linear
structure then the linearity of $D\pi$ is easily proved from the
definitions; hence $\pi$ is smooth.   \qed
\enddemo

\section 3. The de Rham cohomology

A {\it vector field} on $X$ is a smooth section of the tangent bundle
$\Tn^1 X$. The set of vector fields is denoted~$\Gamma(X)$. Given a
vector field $\xi$ and $f$ in $\cader(X)$, define
$$
\xi(f)F := \frac{d}{dt} f(\xi(F)t) \biggr|_{t=0}.
$$

\proclaim Lemma 3.1.
$\xi(f)$ is smooth and therefore $\xi$ is a derivation on $\cader(X)$.

\noindent{\it Proof}.\enspace
The function $\xi(f)$ is well defined since $f$ is smooth. Let
$p \in P(X)$, then $\xi \circ p$ is a plaque on $\Tn^1 X$, that is,
$(\xi \circ p)r = [\bar p(r,t)]_t$. It follows that
$$
\xi(f) \circ p(r) = \xi(p(r))f
 = \frac{d}{dt} f\circ p(r,t) \biggr|_{t=0}.
$$
This is a smooth function since $f$ is differentiable, therefore
$\xi(f)$ is differentiable.

Assume $p_1 \sim^n_F p_2$, then since $\xi$ is smooth we have that
$(\xi \circ p_1)r \approx (\xi \circ p_2)r$ at $\xi(F)$. Let
$\bar p_1$, $\bar p_2$ be such that
$[\bar p_i(r,t)]_t = (\xi \circ p_i)(r)$; then
$\bar p_1(r,t) \sim \bar p_2(r,t)$ at~$F$.

Then
$$
D^n \frac{d}{dt} f(\bar p_1(r,t)) \biggr|_{t=0}
= D^n \frac{d}{dt} f(\bar p_2(r,t))  \biggr|_{t=0}
$$
and so $\xi(f)p_1 \sim^n \xi(f)p_2$. If $X$ is linear, the linearity
of $D^m\xi (f)$ follows from the linearity of $D^m\xi$.   \qed
\enddemo

It is easy to prove that if $X$ has a continuous linear structure,
then $\Gamma(X)$ is a $\cader(X)$-module.

A {\it local flow\/} is a map $\phi \: P(X) \to P(X)$ such that
(a)~if $p \in P^n(X)$ then $\phi(p) \in P^{n+1}(X)$;
(b)~$\phi(p)(r,0) = p(r)$; and (c)~if $p_1(r) = p_2(s)$ then
$[\phi(p_1)(r,t)]_t = [\phi(p_2)(s,t)]_t$. Each local flow defines
a unique vector field and each vector field defines a local flow
$\phi$. Indeed, for each vector field $\xi$ define
$\phi_\xi(p) := \bar p$ (since $\xi$ is smooth) where
$\xi(p(r)) = [\bar p(r,t)]_t$; and conversely, given $\phi$, define
$\xi_\phi(F) = [\phi(p)(r_0,t)]_t$ where $p$ is any plaque such that
$p(r_0) = F$.

Let us assume that $X$ has a continuous linear structure. Let
$\xi_1,\xi_2 \in \Gamma(X)$. Then $B(\xi_1,\xi_2)$, defined by
$$
B(\xi_1,\xi_2)(f) := \xi_1\xi_2 f - \xi_2\xi_1 f,
$$
is an element of $Der(\cader(X))$. Now, let $\M$ be a maximal
subalgebra of $\Gamma(X)$ with this product. For each such subalgebra,
define an $n$-form as a section $\omega$ of $\W^n(\Tn(X))$ such that
$\omega(\xi_1,\dots,\xi_n)(F):= \omega(F)(\xi_1(F),\dots,\xi_n(F))$ is
smooth whenever $\xi_i \in \M$ for all $i = 1,\dots,n$. The set of
these, denoted by $\W^n(X)$, is a $\cader(X)$-module and an
associative algebra with the operation:
$$
\omega \w \eta := \frac{(k +l)!}{k!l!} \Alt(\omega \ox \eta).
$$

An example of a subspace of $\W^n(X)$ is the set of forms $\omega$
expressible as
$$
\omega = \sum_{I = \{i_1,\dots,i_n\}}
  h_I \,df_{i_1} \w\cdots\w df_{i_n},
$$
where the $f_i$ are smooth and $\{h_I\}$ is a locally finite family
of smooth functions.

Define $d_n \: \W^n(X) \to \W^{n+1}(X)$ by
$$
\eqalign{
d_n\omega(\xi_1,\dots,&\xi_{n+1})
= \sum_{i=1}^{n+1} (-1)^{i+1}
     \xi_i\, \omega(\xi_1,\dots,\widehat\xi_i,\dots,\xi_{n+1})
\cr
&+ \sum_{i<j} (-1)^{i+j}
     \omega(B(\xi_i,\xi_j),\xi_1,\dots,\widehat\xi_i,\dots,
            \widehat\xi_j,\dots,\xi_{n+1}).
\cr}
$$

It is clear that $d_{n+1} \circ d_n = 0$. This allows the definition
of $Z_n(X,\R) := \ker(d_n)$, $B_n(X,\R) := \im(d_{n-1})$ and
$$
H^n_\dR(X,\R) := Z_n/B_n.
$$
This is the $n$-th de~Rham cohomology group of~$X$.


\varsection  Acknowledgments

\medskip
I gratefully acknowledge the support of the Vicerrector\'{\i}a de
Investigaci\'on of the University of Costa Rica. I would like also to
thank Joseph C. V\'arilly for helpful discussions.

\bigskip

\varsection  References

\frenchspacing

\refno \Adams.
M. Adams, T. Ratiu and R. Schmid,
A Lie group structure of diffeomorphism groups and invertible
Fourier integral operators, with applications,
in: V. Kac, ed., {\it Infinite Dimensional Groups with Applications},
MSRI Publications {\bf 4} (Springer, Berlin, 1985) 1--26.

\refno \DonatoGeometrie.
P. Donato,
G\'eom\'etrie des orbites coadjointes des groupes de
dif\-feo\-mor\-phismes,
in: C. Albert, ed., {\it G\'eom\'etrie Symplectique et M\'ecanique},
Lecture Notes in Mathematics {\bf 1416} (Springer, Berlin, 1988)
84--104.

\refno \DonatoVirasoro.
P. Donato,
Les diffeomorphismes du circle comme orbit symplectique dans les
moments de Virasoro,
Preprint CPT--92/P.2681, CNRS--Luminy, 1992.

\refno \DonatoFibrations.
P. Donato and P. Iglesias,
Cohomologie des formes dans les espaces diffeologiques,
Preprint CPT--87/P.1986, CNRS--Luminy, 1987.

\refno \KirillovBook.
A. A. Kirillov,
{\it Elements of the Theory of Representations\/}
(Sprin\-ger, Berlin, 1976).

\refno \SouriauAlgoritme.
J. M. Souriau,
Un algoritme g\'en\'erateur de structures quantiques,
Ast\'erisque, hors s\'erie (1985) 341--399.

\refno \Woodhouse.
N. Woodhouse,
{\it Geometric Quantization\/}
(Clarendon Press, Oxford, 1992).

\bye